\documentclass{conm-p-l}

\theoremstyle{definition}

\theoremstyle{remark}

\numberwithin{equation}{section}

\usepackage{graphicx,amsthm}
\usepackage{verbatim,amsmath,amscd,amssymb}
\input xy
\xyoption{all}
\begin{document}
\renewcommand{\labelenumi}{$($\roman{enumi}$)$}
\renewcommand{\labelenumii}{$(${\rm \alph{enumii}}$)$}
\font\germ=eufm10
\newcommand{\cI}{{\mathcal I}}
\newcommand{\cA}{{\mathcal A}}
\newcommand{\cB}{{\mathcal B}}
\newcommand{\cC}{{\mathcal C}}
\newcommand{\cD}{{\mathcal D}}
\newcommand{\cE}{{\mathcal E}}
\newcommand{\cF}{{\mathcal F}}
\newcommand{\cG}{{\mathcal G}}
\newcommand{\cH}{{\mathcal H}}
\newcommand{\cK}{{\mathcal K}}
\newcommand{\cL}{{\mathcal L}}
\newcommand{\cM}{{\mathcal M}}
\newcommand{\cN}{{\mathcal N}}
\newcommand{\cO}{{\mathcal O}}
\newcommand{\cR}{{\mathcal R}}
\newcommand{\cS}{{\mathcal S}}
\newcommand{\cV}{{\mathcal V}}
\newcommand{\fra}{\mathfrak a}
\newcommand{\frb}{\mathfrak b}
\newcommand{\frc}{\mathfrak c}
\newcommand{\frd}{\mathfrak d}
\newcommand{\fre}{\mathfrak e}
\newcommand{\frf}{\mathfrak f}
\newcommand{\frg}{\mathfrak g}
\newcommand{\frh}{\mathfrak h}
\newcommand{\fri}{\mathfrak i}
\newcommand{\frj}{\mathfrak j}
\newcommand{\frk}{\mathfrak k}
\newcommand{\frI}{\mathfrak I}
\newcommand{\fm}{\mathfrak m}
\newcommand{\frn}{\mathfrak n}
\newcommand{\frp}{\mathfrak p}
\newcommand{\fq}{\mathfrak q}
\newcommand{\frr}{\mathfrak r}
\newcommand{\frs}{\mathfrak s}
\newcommand{\frt}{\mathfrak t}
\newcommand{\fru}{\mathfrak u}
\newcommand{\frA}{\mathfrak A}
\newcommand{\frB}{\mathfrak B}
\newcommand{\frF}{\mathfrak F}
\newcommand{\frG}{\mathfrak G}
\newcommand{\frH}{\mathfrak H}
\newcommand{\frJ}{\mathfrak J}
\newcommand{\frN}{\mathfrak N}
\newcommand{\frP}{\mathfrak P}
\newcommand{\frT}{\mathfrak T}
\newcommand{\frU}{\mathfrak U}
\newcommand{\frV}{\mathfrak V}
\newcommand{\frX}{\mathfrak X}
\newcommand{\frY}{\mathfrak Y}
\newcommand{\frZ}{\mathfrak Z}
\newcommand{\rA}{\mathrm{A}}
\newcommand{\rC}{\mathrm{C}}
\newcommand{\rd}{\mathrm{d}}
\newcommand{\rB}{\mathrm{B}}
\newcommand{\rD}{\mathrm{D}}
\newcommand{\rE}{\mathrm{E}}
\newcommand{\rH}{\mathrm{H}}
\newcommand{\rK}{\mathrm{K}}
\newcommand{\rL}{\mathrm{L}}
\newcommand{\rM}{\mathrm{M}}
\newcommand{\rN}{\mathrm{N}}
\newcommand{\rR}{\mathrm{R}}
\newcommand{\rT}{\mathrm{T}}
\newcommand{\rZ}{\mathrm{Z}}
\newcommand{\bbA}{\mathbb A}
\newcommand{\bbB}{\mathbb B}
\newcommand{\bbC}{\mathbb C}
\newcommand{\bbG}{\mathbb G}
\newcommand{\bbF}{\mathbb F}
\newcommand{\bbH}{\mathbb H}
\newcommand{\bbP}{\mathbb P}
\newcommand{\bbN}{\mathbb N}
\newcommand{\bbQ}{\mathbb Q}
\newcommand{\bbR}{\mathbb R}
\newcommand{\bbV}{\mathbb V}
\newcommand{\bbZ}{\mathbb Z}
\newcommand{\adj}{\operatorname{adj}}
\newcommand{\Ad}{\mathrm{Ad}}
\newcommand{\Ann}{\mathrm{Ann}}
\newcommand{\rcris}{\mathrm{cris}}
\newcommand{\ch}{\mathrm{ch}}
\newcommand{\coker}{\mathrm{coker}}
\newcommand{\diag}{\mathrm{diag}}
\newcommand{\Diff}{\mathrm{Diff}}
\newcommand{\Dist}{\mathrm{Dist}}
\newcommand{\rDR}{\mathrm{DR}}
\newcommand{\ev}{\mathrm{ev}}
\newcommand{\Ext}{\mathrm{Ext}}
\newcommand{\cExt}{\mathcal{E}xt}
\newcommand{\fin}{\mathrm{fin}}
\newcommand{\Frac}{\mathrm{Frac}}
\newcommand{\GL}{\mathrm{GL}}
\newcommand{\Hom}{\mathrm{Hom}}
\newcommand{\hd}{\mathrm{hd}}
\newcommand{\rht}{\mathrm{ht}}
\newcommand{\id}{\mathrm{id}}
\newcommand{\im}{\mathrm{im}}
\newcommand{\inc}{\mathrm{inc}}
\newcommand{\ind}{\mathrm{ind}}
\newcommand{\coind}{\mathrm{coind}}
\newcommand{\Lie}{\mathrm{Lie}}
\newcommand{\Max}{\mathrm{Max}}
\newcommand{\mult}{\mathrm{mult}}
\newcommand{\op}{\mathrm{op}}
\newcommand{\ord}{\mathrm{ord}}
\newcommand{\pt}{\mathrm{pt}}
\newcommand{\qt}{\mathrm{qt}}
\newcommand{\rad}{\mathrm{rad}}
\newcommand{\res}{\mathrm{res}}
\newcommand{\rgt}{\mathrm{rgt}}
\newcommand{\rk}{\mathrm{rk}}
\newcommand{\SL}{\mathrm{SL}}
\newcommand{\soc}{\mathrm{soc}}
\newcommand{\Spec}{\mathrm{Spec}}
\newcommand{\St}{\mathrm{St}}
\newcommand{\supp}{\mathrm{supp}}
\newcommand{\Tor}{\mathrm{Tor}}
\newcommand{\Tr}{\mathrm{Tr}}
\newcommand{\wt}{\mathrm{wt}}
\newcommand{\Ab}{\mathbf{Ab}}
\newcommand{\Alg}{\mathbf{Alg}}
\newcommand{\Grp}{\mathbf{Grp}}
\newcommand{\Mod}{\mathbf{Mod}}
\newcommand{\Sch}{\mathbf{Sch}}\newcommand{\bfmod}{{\bf mod}}
\newcommand{\Qc}{\mathbf{Qc}}
\newcommand{\Rng}{\mathbf{Rng}}
\newcommand{\Top}{\mathbf{Top}}
\newcommand{\Var}{\mathbf{Var}}
\newcommand{\gromega}{\langle\omega\rangle}
\newcommand{\lbr}{\begin{bmatrix}}
\newcommand{\rbr}{\end{bmatrix}}
\newcommand{\forb}{\bigcirc\kern-2.8ex \because}
\newcommand{\forbb}{\bigcirc\kern-3.0ex \because}
\newcommand{\forbbb}{\bigcirc\kern-3.1ex \because}
\newcommand{\cd}{commutative diagram }
\newcommand{\SpS}{spectral sequence}
\newcommand\C{\mathbb C}
\newcommand\hh{{\hat{H}}}
\newcommand\eh{{\hat{E}}}
\newcommand\F{\mathbb F}
\newcommand\fh{{\hat{F}}}
\def\ge{\frg}
\def\AA{{\mathcal A}}
\def\al{\alpha}
\def\bq{B_q(\ge)}
\def\bqm{B_q^-(\ge)}
\def\bqz{B_q^0(\ge)}
\def\bqp{B_q^+(\ge)}
\def\beneme{\begin{enumerate}}
\def\beq{\begin{equation}}
\def\beqn{\begin{eqnarray}}
\def\beqnn{\begin{eqnarray*}}
\def\bigsl{{\hbox{\fontD \char'54}}}
\def\bbra#1,#2,#3{\left\{\begin{array}{c}\hspace{-5pt}
#1;#2\\ \hspace{-5pt}#3\end{array}\hspace{-5pt}\right\}}
\def\cd{\cdots}
\def\CC{\mathbb{C}}
\def\CBL{\cB_L(\TY(B,1,n+1))}
\def\CBM{\cB_M(\TY(B,1,n+1))}
\def\CVL{\cV_L(\TY(D,1,n+1))}
\def\CVM{\cV_M(\TY(D,1,n+1))}
\def\ddd{\hbox{\germ D}}
\def\del{\delta}
\def\Del{\Delta}
\def\Delr{\Delta^{(r)}}
\def\Dell{\Delta^{(l)}}
\def\Delb{\Delta^{(b)}}
\def\Deli{\Delta^{(i)}}
\def\Delre{\Delta^{\rm re}}
\def\ei{e_i}
\def\eit{\tilde{e}_i}
\def\eneme{\end{enumerate}}
\def\ep{\epsilon}
\def\eeq{\end{equation}}
\def\eeqn{\end{eqnarray}}
\def\eeqnn{\end{eqnarray*}}
\def\fit{\tilde{f}_i}
\def\FF{{\rm F}}
\def\ft{\tilde{f}}
\def\gau#1,#2{\left[\begin{array}{c}\hspace{-5pt}#1\\
\hspace{-5pt}#2\end{array}\hspace{-5pt}\right]}
\def\gl{\hbox{\germ gl}}
\def\hom{{\hbox{Hom}}}
\def\ify{\infty}
\def\io{\iota}
\def\kp{k^{(+)}}
\def\km{k^{(-)}}
\def\llra{\relbar\joinrel\relbar\joinrel\relbar\joinrel\rightarrow}
\def\lan{\langle}
\def\lar{\longrightarrow}
\def\max{{\rm max}}
\def\lm{\lambda}
\def\Lm{\Lambda}
\def\mapright#1{\smash{\mathop{\longrightarrow}\limits^{#1}}}
\def\Mapright#1{\smash{\mathop{\Longrightarrow}\limits^{#1}}}
\def\mm{{\bf{\rm m}}}
\def\nd{\noindent}
\def\nn{\nonumber}
\def\nnn{\hbox{\germ n}}
\def\catob{{\mathcal O}(B)}
\def\oint{{\mathcal O}_{\rm int}(\ge)}
\def\ot{\otimes}
\def\op{\oplus}
\def\opi{\ovl\pi_{\lm}}
\def\osigma{\ovl\sigma}
\def\ovl{\overline}
\def\plm{\Psi^{(\lm)}_{\io}}
\def\qq{\qquad}
\def\q{\quad}
\def\qed{\hfill\framebox[2mm]{}}
\def\QQ{\mathbb Q}
\def\qi{q_i}
\def\qii{q_i^{-1}}
\def\ra{\rightarrow}
\def\ran{\rangle}
\def\rlm{r_{\lm}}
\def\ssl{\hbox{\germ sl}}
\def\slh{\widehat{\ssl_2}}
\def\ti{t_i}
\def\tii{t_i^{-1}}
\def\til{\tilde}
\def\tm{\times}
\def\tt{\frt}
\def\TY(#1,#2,#3){#1^{(#2)}_{#3}}
\def\ua{U_{\AA}}
\def\ue{U_{\vep}}
\def\uq{U_q(\ge)}
\def\uqp{U'_q(\ge)}
\def\ufin{U^{\rm fin}_{\vep}}
\def\ufinp{(U^{\rm fin}_{\vep})^+}
\def\ufinm{(U^{\rm fin}_{\vep})^-}
\def\ufinz{(U^{\rm fin}_{\vep})^0}
\def\uqm{U^-_q(\ge)}
\def\uqmq{{U^-_q(\ge)}_{\bf Q}}
\def\uqpm{U^{\pm}_q(\ge)}
\def\uqq{U_{\bf Q}^-(\ge)}
\def\uqz{U^-_{\bf Z}(\ge)}
\def\ures{U^{\rm res}_{\AA}}
\def\urese{U^{\rm res}_{\vep}}
\def\uresez{U^{\rm res}_{\vep,\ZZ}}
\def\util{\widetilde\uq}
\def\uup{U^{\geq}}
\def\ulow{U^{\leq}}
\def\bup{B^{\geq}}
\def\blow{\ovl B^{\leq}}
\def\vep{\varepsilon}
\def\vp{\varphi}
\def\vpi{\varphi^{-1}}
\def\VV{{\mathcal V}}
\def\xii{\xi^{(i)}}
\def\Xiioi{\Xi_{\io}^{(i)}}
\def\W1{W(\varpi_1)}
\def\WW{{\mathcal W}}
\def\wt{{\rm wt}}
\def\wtil{\widetilde}
\def\what{\widehat}
\def\wpi{\widehat\pi_{\lm}}
\def\ZZ{\mathbb Z}

\def\m@th{\mathsurround=0pt}
\def\fsquare(#1,#2){
\hbox{\vrule$\hskip-0.4pt\vcenter to #1{\normalbaselines\m@th
\hrule\vfil\hbox to #1{\hfill$\scriptstyle #2$\hfill}\vfil\hrule}$\hskip-0.4pt
\vrule}}

\newtheorem{thm}{Theorem}[section]
\newtheorem{pro}[thm]{Proposition}
\newtheorem{lem}[thm]{Lemma}
\newtheorem{ex}[thm]{Example}
\newtheorem{cor}[thm]{Corollary}
\newtheorem{conj}[thm]{Conjecture}
\theoremstyle{definition}
\newtheorem{df}[thm]{Definition}

\newcommand{\cmt}{\marginpar}
\newcommand{\seteq}{\mathbin{:=}}
\newcommand{\cl}{\colon}
\newcommand{\be}{\begin{enumerate}}
\newcommand{\ee}{\end{enumerate}}
\newcommand{\bnum}{\be[{\rm (i)}]}
\newcommand{\enum}{\ee}
\newcommand{\ro}{{\rm(}}
\newcommand{\rf}{{\rm)}}
\newcommand{\set}[2]{\left\{#1\,\vert\,#2\right\}}
\newcommand{\sbigoplus}{{\mbox{\small{$\bigoplus$}}}}
\newcommand{\ba}{\begin{array}}
\newcommand{\ea}{\end{array}}
\newcommand{\on}{\operatorname}
\newcommand{\eq}{\begin{eqnarray}}
\newcommand{\eneq}{\end{eqnarray}}
\newcommand{\hs}{\hspace*}

\title{Affine Geometric Crystal of type $\TY(D,3,4)$}

\author[M.Igarashi]{Mana I\textsc{garashi}}
\author[T.Nakashima]{Toshiki N\textsc{akashima}*}
\address{Department of Mathematics, 
Sophia University, Kioicho 7-1, Chiyoda-ku, Tokyo 102-8554,
Japan}
\email{mana-i@hoffman.cc.sophia.ac.jp}
\address{Department of Mathematics, 
Sophia University, Kioicho 7-1, Chiyoda-ku, Tokyo 102-8554,
Japan}
\email{toshiki@mm.sophia.ac.jp}
\thanks{*supported in part by JSPS Grants 
in Aid for Scientific Research \#19540050.}
\subjclass{Primary 17B37; 17B67; 
Secondary 22E65; 14M15}
\date{}


\keywords{affine geometric crystal, 
positive structure,type $\TY(D,3,4)$}

\begin{abstract}
We shall realize certain affine geometric 
crystal of type $\TY(D,3,4)$ 
associated with the fundamental 
representation $W(\varpi_1)$ explicitly . 
By its explicit form, we see that it 
has a positive structure.
\end{abstract}

\maketitle
\renewcommand{\thesection}{\arabic{section}}
\section{Introduction}
\setcounter{equation}{0}
\renewcommand{\theequation}{\thesection.\arabic{equation}}

The notion of geometric crystals is introduced as 
a geometric analogue to Kashiwara's crystals(\cite{BK}).
For a fixed Cartan data $(A,\{\al_i\}_{i\in I},\{h_i\}_{\i\in I})$, it
is defined as a quadruple $(X,\{e_i\}_{i\in I},
\{\gamma_i\}_{i\in I}, \{\vep_i\}_{i\in I})$
where X is 
an algebraic(ind-)variety over the complex number $\bbC$, 
$e_i$ is a rational $\bbC^\times$-action 
$e_i:\bbC^\times\times X\longrightarrow X$ and 
 $\gamma_i,\vep_i
:X\longrightarrow \bbC$ $(i\in I)$ are 
rational functions satisfying certain conditions 
(see Definition \ref{def-gc}).
One of the remarkable properties of 
geometric crystal is that 
if they are equipped with the so-called 
``positive structure'',
there exists a functor from 
certain category of geometric crystals to
the category of Langlands dual crystals, say, 
tropicalization/
ultra-discretization procedure(see \ref{positive-str}).

In \cite{KNO}, we gave  conjectures for  
constructions of  some affine geometric crystals and 
their relations to limit of perfect crystals.
Therein, some partial answers are presented by 
explicitly constructing affine geometric crystals. 
Adopting the same method, we obtained 
the affine geometric crystal $\cV$ of type $\TY(G,1,2)$
(\cite{N3}) very explicitly and 
see that it has a positive structure.
In \cite{N4}, it is shown that 
its ultra-discretization is isomorphic to 
certain limit of perfect crystals of type 
$\TY(D,3,4)$ (\cite{KMOY}).

In this article, we construct an affine 
geometric crystal of type $\TY(D,3,4)$ associated with 
the fundamental representation $W(\varpi_1)$ by
the same way as the $\TY(G,1,2)$-case in \cite{N3}
and see its positive structure.
Though we have the positive structure, 
in this article 
we do not treat its ultra-discertization to the 
corresponding $\TY(G,1,2)$-crystals,
which will be discussed elsewhere.

Let us explain how to obtain the 
affine geometric crystal $\cV$ of type $\TY(D,3,4)$:
Let $I:=\{0,1,2\}$ be the index set of simple roots and
$\{\Lm_i\}_{i\in I}$ the set of fundamental weights
(see \ref{KM}). 
Let $\varpi_1:=\Lm_1-2\Lm_0$ be 
the first level 0 fundamental weight
and $W(\varpi_1)$ the associated 
fundamental representation (see \ref{fundamental}), which
is an 8-dimensional module with the basis
$M=\{v_1,v_2,v_3,v_0,\emptyset, v_{\ovl 3},
v_{\ovl 2},v_{\ovl 1}\}$. Set 
\begin{eqnarray*}
&&\hspace{-30pt}\cV_1:=\{V_1(x)
:=Y_0(x_0)Y_1(x_1)Y_2(x_2)Y_1(x_3)Y_2(x_4)Y_1(x_5)
v_1\,\,\vert\,\,x_i\in\bbC^\times,(0\leq i\leq 5)\},\\
&&\hspace{-30pt}\cV_2:=\{V_2(y):=
Y_2(y_2)Y_1(y_1)Y_2(y_4)Y_1(y_3)Y_0(y_0)Y_1(y_5)
v_{\ovl 2}\,\,\vert\,\,y_i\in\bbC^\times,
(0\leq i\leq 5)\},
\end{eqnarray*}
where $Y_i(c):=y_i(c^{-1})\al^\vee_i(c)$ (see \ref{KM})
and note that $\cV_1,\cV_2\subset W(\varpi_1)$ and 
$\cV_1$ (resp. $\cV_2$) has a 
$(\TY(D,3,4))_{\{1,2\}}(\cong G_2)$ 
(resp. $(\TY(D,3,4))_{\{0,1\}}(\cong A_2)$)-
geometric crystal structure.
Each $V_k\in \cV_k$ $(k=1,2)$ is in the form:
\[
 V_1(x):=\sum_{m\in M}X_m m,\q
 V_2(y):=\sum_{m\in M}Y_m m,
\]
where $X_m$ (resp. $Y_m$) is a rational function
in $(x_0,\cd,x_5)$ (resp. $(y_0,\cd,y_5)$).
First, for given $x$ we solve the equation 
\[
 V_2(y)=a(x)V_1(x),
\]
where $a(x)$ is a rational function.
Then we obtain the unique solution 
$y=\osigma(x)$ and $a(x)$. This $\osigma$ defines a
rational map  from $\cV_1$ to $\cV_2$ 
($V_1(x)\mapsto V_2(y):=\osigma(V_1(x))$).
Next, we shall see that this rational map is 
bi-positive and birational.
Then, we define $e^c_0,\gamma_0,\vep_0$ on 
$\cV_1$ by $e^c_0(x):=
\osigma^{-1}\circ e^c_0\circ\osigma(x)$,
$\gamma_0(x):=\gamma_0(\osigma(x))$ and 
$\vep_0(x):=\vep_0(\osigma(x))$.
This gives a $\TY(D,3,4)$-geometric 
crystal structure on $\cV_1$.
Finally, we present Conjecture \ref{conjecture}
as a further problem.

\renewcommand{\thesection}{\arabic{section}}
\section{Geometric crystals}
\setcounter{equation}{0}
\renewcommand{\theequation}{\thesection.\arabic{equation}}

In this section, 
we review Kac-Moody groups and geometric crystals
following 
\cite{PK}, \cite{Ku2}, \cite{BK}
\subsection{Preliminaries and Notations}
\label{KM}
Fix a symmetrizable generalized Cartan matrix
 $A=(a_{ij})_{i,j\in I}$ with a finite index set $I$.
Let $(\tt,\{\al_i\}_{i\in I},\{\al^\vee_i\}_{i\in I})$ 
be the associated
root data, where ${\tt}$ is a vector space 
over $\bbC$ and
$\{\al_i\}_{i\in I}\subset\tt^*$ and 
$\{\al^\vee_i\}_{i\in I}\subset\tt$
are linearly independent 
satisfying $\al_j(\al^\vee_i)=a_{ij}$.

The Kac-Moody Lie algebra $\ge=\ge(A)$ associated with $A$
is the Lie algebra over $\bbC$ generated by $\tt$, the 
Chevalley generators $e_i$ and $f_i$ $(i\in I)$
with the usual defining relations (\cite{Kac}).
There is the root space decomposition 
$\ge=\bigoplus_{\al\in \tt^*}\ge_{\al}$.
Denote the set of roots by 
$\Delta:=\{\al\in \tt^*|\al\ne0,\,\,\ge_{\al}\ne(0)\}$.
Set $Q=\sum_i\bbZ \al_i$, $Q_+=\sum_i\bbZ_{\geq0} \al_i$,
$Q^\vee:=\sum_i\bbZ \al^\vee_i$
and $\Delta_+:=\Delta\cap Q_+$.
An element of $\Delta_+$ is called 
a {\it positive root}.
Let $P\subset \tt^*$ be a weight lattice such that 
$\bbC\ot P=\tt^*$, whose element is called a
weight.

Define simple reflections $s_i\in{\rm Aut}(\tt)$ $(i\in I)$ by
$s_i(h):=h-\al_i(h)\al^\vee_i$, which generate the Weyl group $W$.
It induces the action of $W$ on $\tt^*$ by
$s_i(\lm):=\lm-\lm(\al^\vee_i)\al_i$.
Set $\Delre:=\{w(\al_i)|w\in W,\,\,i\in I\}$, whose element 
is called a real root.

Let $\ge'$ be the derived Lie algebra 
of $\ge$ and let 
$G$ be the Kac-Moody group associated 
with $\ge'$(\cite{PK}).
Let $U_{\al}:=\exp\ge_{\al}$ $(\al\in \Delre)$
be the one-parameter subgroup of $G$.
The group $G$ is generated by $U_{\al}$ $(\al\in \Delre)$.
Let $U^{\pm}$ be the subgroup generated by $U_{\pm\al}$
($\al\in \Delre_+=\Delre\cap Q_+$), {\it i.e.,}
$U^{\pm}:=\lan U_{\pm\al}|\al\in\Del^{\rm re}_+\ran$.

For any $i\in I$, there exists a unique homomorphism;
$\phi_i:SL_2(\bbC)\rightarrow G$ such that
\[
\hspace{-2pt}\phi_i\left(
\left(
\begin{array}{cc}
c&0\\
0&c^{-1}
\end{array}
\right)\right)=c^{\al^\vee_i},\,
\phi_i\left(
\left(
\begin{array}{cc}
1&t\\
0&1
\end{array}
\right)\right)=\exp(t e_i),\,
 \phi_i\left(
\left(
\begin{array}{cc}
1&0\\
t&1
\end{array}
\right)\right)=\exp(t f_i).
\]
where $c\in\bbC^\times$, $t\in\bbC$ and 
$c^{\al^\vee_i}\in \operatorname{Hom}_{\bbZ}(P,\bbC^\times)$ such that 
$c^{\al^\vee_i}(\lm)=c^{\lm(\al^\vee_i)}$.
Set $\al^\vee_i(c):=c^{\al^\vee_i}$,
$x_i(t):=\exp{(t e_i)}$, $y_i(t):=\exp{(t f_i)}$, 
$G_i:=\phi_i(SL_2(\bbC))$,
$T_i:=\phi_i(\{{\rm diag}(c,c^{-1})\vert 
c\in\bbC^{\vee}\})$ 
and 
$N_i:=N_{G_i}(T_i)$. Let
$T$ (resp. $N$) be the subgroup of $G$ 
with the Lie algebra $\tt$
(resp. generated by the $N_i$'s), 
which is called a {\it maximal torus} in $G$, and let
$B^{\pm}=U^{\pm}T$ be the Borel subgroup of $G$.
We have the isomorphism
$\phi:W\mapright{\sim}N/T$ defined by $\phi(s_i)=N_iT/T$.
An element $\ovl s_i:=x_i(-1)y_i(1)x_i(-1)$ is in 
$N_G(T)$, which is a representative of 
$s_i\in W=N_G(T)/T$. 

\subsection{Geometric crystals}

Let $X$ be an algebraic(ind)-variety , 
{$\gamma_i:X\rightarrow \bbC$} and 
$\vep_i:X\longrightarrow \bbC$ ($i\in I$) 
rational functions on $X$, and
{$e_i:\bbC^\times \times X\longrightarrow X$}
$((c,x)\mapsto e^c_i(x))$ a
rational $\bbC^\times$-action.

\begin{df}
\label{def-gc}
A quadruple $(X,\{e_i\}_{i\in I},\{\gamma_i,\}_{i\in I},
\{\vep_i\}_{i\in I})$ is a 
$G$ (or $\ge$)-\\{\it geometric} {\it crystal} 
if
\begin{enumerate}
\item
$\{1\}\times X\subset dom(e_i)$ 
for any $i\in I$.
\item
$\gamma_j(e^c_i(x))=c^{a_{ij}}\gamma_j(x)$.
\item
The rational $\bbC^\times$ actions $\{e_i\}_{i\in I}$ satisfy the
following relations ({\it Verma relations}):
\[
 \begin{array}{lll}
&\hspace{-20pt}e^{c_1}_{i}e^{c_2}_{j}
=e^{c_2}_{j}e^{c_1}_{i}&
{\rm if }\,\,a_{ij}=a_{ji}=0,\\
&\hspace{-20pt} e^{c_1}_{i}e^{c_1c_2}_{j}e^{c_2}_{i}
=e^{c_2}_{j}e^{c_1c_2}_{i}e^{c_1}_{j}&
{\rm if }\,\,a_{ij}=a_{ji}=-1,\\
&\hspace{-20pt}
e^{c_1}_{i}e^{c^2_1c_2}_{j}e^{c_1c_2}_{i}e^{c_2}_{j}
=e^{c_2}_{j}e^{c_1c_2}_{i}e^{c^2_1c_2}_{j}e^{c_1}_{i}&
{\rm if }\,\,a_{ij}=-2,\,
a_{ji}=-1,\\
&\hspace{-20pt}
e^{c_1}_{i}e^{c^3_1c_2}_{j}e^{c^2_1c_2}_{i}
e^{c^3_1c^2_2}_{j}e^{c_1c_2}_{i}e^{c_2}_{j}
=e^{c_2}_{j}e^{c_1c_2}_{i}e^{c^3_1c^2_2}_{j}e^{c^2_1c_2}_{i}
e^{c^3_1c_2}_je^{c_1}_i&
{\rm if }\,\,a_{ij}=-3,\,
a_{ji}=-1,
\end{array}
\]
\item
$\vep_i(e_i^c(x))=c^{-1}\vep_i(x)$.
\end{enumerate}
\end{df}
Note that the last formula in (iii) is different from the one in 
\cite{BK}, \cite{N}, \cite{N2} which seems to be
incorrect.

\subsection{Geometric crystal on Schubert cell}
\label{schubert}

Let $w\in W$ be a Weyl group element and take a 
reduced expression $w=s_{i_1}\cd s_{i_l}$. 
Let $X:=G/B$ be the flag
variety, which is an ind-variety 
and $X_w\subset X$ the
Schubert cell associated with $w$, which has 
a natural geometric crystal structure
(\cite{BK},\cite{N}).
For ${\bf i}:=(i_1,\cd,i_k)$, set 
\begin{equation}
B_{\bf i}^-
:=\{Y_{\bf i}(c_1,\cd,c_k)
:=Y_{i_1}(c_1)\cd Y_{i_l}(c_k)
\,\vert\, c_1\cd,c_k\in\bbC^\times\}\subset B^-,
\label{bw1}
\end{equation}
which has a geometric crystal structure(\cite{N})
isomorphic to $X_w$. 
The explicit forms of the action $e^c_i$, the rational 
function $\vep_i$  and $\gamma_i$ on 
$B_{\bf i}^-$ are given by
\begin{eqnarray}
&& e_i^c(Y_{i_1}(c_1)\cd Y_{i_l}(c_k))
=Y_{i_1}({\mathcal C}_1)\cd Y_{i_l}({\mathcal C}_k)),\nn \\
&&\text{where}\nn\\
&&{\mathcal C}_j:=
c_j\cdot \frac{\displaystyle \sum_{1\leq m\leq j,i_m=i}
 \frac{c}
{c_1^{a_{i_1,i}}\cd c_{m-1}^{a_{i_{m-1},i}}c_m}
+\sum_{j< m\leq k,i_m=i} \frac{1}
{c_1^{a_{i_1,i}}\cd c_{m-1}^{a_{i_{m-1},i}}c_m}}
{\displaystyle\sum_{1\leq m<j,i_m=i} 
 \frac{c}
{c_1^{a_{i_1,i}}\cd c_{m-1}^{a_{i_{m-1},i}}c_m}+
\mathop\sum_{j\leq m\leq k,i_m=i}  \frac{1}
{c_1^{a_{i_1,i}}\cd c_{m-1}^{a_{i_{m-1},i}}c_m}},
\label{eici}\\
&& \vep_i(Y_{i_1}(c_1)\cd Y_{i_l}(c_k))=
\sum_{1\leq m\leq k,i_m=i} \frac{1}
{c_1^{a_{i_1,i}}\cd c_{m-1}^{a_{i_{m-1},i}}c_m},
\label{vep-i}\\
&&\gamma_i(Y_{i_1}(c_1)\cd Y_{i_l}(c_k))
=c_1^{a_{i_1,i}}\cd c_k^{a_{i_k,i}}.
\label{gamma-i}
\end{eqnarray}

\subsection{Positive structure,\,\,
Ultra-discretizations \,\, and \,\,Tropicalizations}
\label{positive-str}
Let us recall the notions of 
positive structure, ultra-discretization and tropicalization.

The setting below is same as \cite{KNO}.
Let $T=(\bbC^\times)^l$ be an algebraic torus over $\bbC$ and 
$X^*(T):={\rm Hom}(T,\bbC^\times)\cong \ZZ^l$ 
(resp. $X_*(T):={\rm Hom}(\bbC^\times,T)\cong \ZZ^l$) 
be the lattice of characters
(resp. co-characters)
of $T$. 
Set $R:=\bbC(c)$ and define
$$
\begin{array}{cccc}
v:&R\setminus\{0\}&\longrightarrow &\ZZ\\
&f(c)&\mapsto
&{\rm deg}(f(c)),
\end{array}
$$
where $\rm deg$ is the degree of poles at $c=\ify$. 
Here note that for $f_1,f_2\in R\setminus\{0\}$, we have
\begin{equation}
v(f_1 f_2)=v(f_1)+v(f_2),\q
v\left(\frac{f_1}{f_2}\right)=v(f_1)-v(f_2)
\label{ff=f+f}
\end{equation}
A non-zero rational function on
an algebraic torus $T$ is called {\em positive} if
it is written as $g/h$ where
$g$ and $h$ are a positive linear combination of
characters of $T$.
\begin{df}
Let 
$f\cl T\rightarrow T'$ be 
a rational morphism between
two algebraic tori $T$ and 
$T'$.
We say that $f$ is {\em positive},
if $\chi\circ f$ is positive
for any character $\chi\cl T'\to \C$.
\end{df}
Denote by ${\rm Mor}^+(T,T')$ the set of 
positive rational morphisms from $T$ to $T'$.

\begin{lem}[\cite{BK}]
\label{TTT}
For any $f\in {\rm Mor}^+(T_1,T_2)$             
and $g\in {\rm Mor}^+(T_2,T_3)$, 
the composition $g\circ f$
is well-defined and belongs to ${\rm Mor}^+(T_1,T_3)$.
\end{lem}

By Lemma \ref{TTT}, we can define a category ${\mathcal T}_+$
whose objects are algebraic tori over $\bbC$ and arrows
are positive rational morphisms.

Let $f\cl T\rightarrow T'$ be a 
positive rational morphism
of algebraic tori $T$ and 
$T'$.
We define a map $\what f\cl X_*(T)\rightarrow X_*(T')$ by 
\[
\langle\chi,\what f(\xi)\rangle
=v(\chi\circ f\circ \xi),
\]
where $\chi\in X^*(T')$ and $\xi\in X_*(T)$.
\begin{lem}[\cite{BK}]
For any algebraic tori $T_1$, $T_2$, $T_3$, 
and positive rational morphisms 
$f\in {\rm Mor}^+(T_1,T_2)$, 
$g\in {\rm Mor}^+(T_2,T_3)$, we have
$\what{g\circ f}=\what g\circ\what f.$
\end{lem}
By this lemma, we obtain a functor 
\[
\begin{array}{cccc}
{\mathcal UD}:&{\mathcal T}_+&\longrightarrow &{{\hbox{\germ Set}}}\\
&T&\mapsto& X_*(T)\\
&(f:T\rightarrow T')&\mapsto& 
(\what f:X_*(T)\rightarrow X_*(T')))
\end{array}
\]


\begin{df}[\cite{BK}]
Let $\chi=(X,\{e_i\}_{i\in I},\{{\rm wt}_i\}_{i\in I},
\{\vep_i\}_{i\in I})$ be a 
geometric crystal, $T'$ an algebraic torus
and $\theta:T'\rightarrow X$ 
a birational isomorphism.
The isomorphism $\theta$ is called 
{\it positive structure} on
$\chi$ if it satisfies
\begin{enumerate}
\item for any $i\in I$ the rational functions
$\gamma_i\circ \theta:T'\rightarrow \bbC$ and 
$\vep_i\circ \theta:T'\rightarrow \bbC$ 
are positive.
\item
For any $i\in I$, the rational morphism 
$e_{i,\theta}:\bbC^\tm \tm T'\rightarrow T'$ defined by
$e_{i,\theta}(c,t)
:=\theta^{-1}\circ e_i^c\circ \theta(t)$
is positive.
\end{enumerate}
\end{df}
Let $\theta:T\rightarrow X$ be a positive structure on 
a geometric crystal $\chi=(X,\{e_i\}_{i\in I},$
$\{{\rm wt}_i\}_{i\in I},
\{\vep_i\}_{i\in I})$.
Applying the functor ${\mathcal UD}$ 
to positive rational morphisms
$e_{i,\theta}:\bbC^\tm \tm T'\rightarrow T'$ and
$\gamma\circ \theta:T'\ra T$
(the notations are
as above), we obtain
\begin{eqnarray*}
\til e_i&:=&{\mathcal UD}(e_{i,\theta}):
\ZZ\tm X_*(T) \rightarrow X_*(T)\\
{\rm wt}_i&:=&{\mathcal UD}(\gamma_i\circ\theta):
X_*(T')\rightarrow \bbZ,\\
\vep_i&:=&{\mathcal UD}(\vep_i\circ\theta):
X_*(T')\rightarrow \bbZ.
\end{eqnarray*}
Now, for given positive structure $\theta:T'\rightarrow X$
on a geometric crystal 
$\chi=(X,\{e_i\}_{i\in I},$
$\{{\rm wt}_i\}_{i\in I},
\{\vep_i\}_{i\in I})$, we associate 
the quadruple $(X_*(T'),\{\til e_i\}_{i\in I},
\{{\rm wt}_i\}_{i\in I},\{\vep_i\}_{i\in I})$
with a free pre-crystal structure (see \cite[2.2]{BK}) 
and denote it by ${\mathcal UD}_{\theta,T'}(\chi)$.
We have the following theorem:

\begin{thm}[\cite{BK}\cite{N}]
For any geometric crystal 
$\chi=(X,\{e_i\}_{i\in I},\{\gamma_i\}_{i\in I},$
$\{\vep_i\}_{i\in I})$ and positive structure
$\theta:T'\rightarrow X$, the associated pre-crystal 
${\mathcal UD}_{\theta,T'}(\chi)=
(X_*(T'),\{e_i\}_{i\in I},\{{\rm wt}_i\}_{i\in I},
\{\vep_i\}_{i\in I})$ 
is a crystal {\rm (see \cite[2.2]{BK})}
\end{thm}

Now, let ${\mathcal GC}^+$ be a category whose 
object is a triplet
$(\chi,T',\theta)$ where 
$\chi=(X,\{e_i\},\{\gamma_i\},\{\vep_i\})$ 
is a geometric crystal and $\theta:T'\rightarrow X$ 
is a positive structure on $\chi$, and morphism
$f:(\chi_1,T'_1,\theta_1)\longrightarrow 
(\chi_2,T'_2,\theta_2)$ is given by a morphism 
$\vp:X_1\longrightarrow X_2$  
($\chi_i=(X_i,\cd)$) such that 
\[
f:=\theta_2^{-1}\circ\vp\circ\theta_1:T'_1\longrightarrow T'_2,
\]
is a positive rational morphism. Let ${\mathcal CR}$
be a category of crystals. 
Then by the theorem above, we have
\begin{cor}
\label{cor-posi}
$\mathcal UD_{\theta,T'}$ as above defines a functor
\begin{eqnarray*}
 {\mathcal UD}&:&{\mathcal GC}^+\longrightarrow {\mathcal CR},\\
&&(\chi,T',\theta)\mapsto X_*(T'),\\
&&(f:(\chi_1,T'_1,\theta_1)\rightarrow 
(\chi_2,T'_2,\theta_2))\mapsto
(\what f:X_*(T'_1)\rightarrow X_*(T'_2)).
\end{eqnarray*}

\end{cor}
We call the functor $\mathcal UD$
{\it ``ultra-discretization''} as \cite{N},\cite{N2}
instead of ``tropicalization'' as in \cite{BK}.
And 
for a crystal $B$, if there
exists a geometric crystal $\chi$ and a positive 
structure $\theta:T'\rightarrow X$ on $\chi$ such that 
${\mathcal UD}(\chi,T',\theta)\cong B$ as crystals, 
we call an object $(\chi,T',\theta)$ in ${\mathcal GC}^+$
a {\it tropicalization} of $B$, where 
it is not known that this correspondence is a functor.

\renewcommand{\thesection}{\arabic{section}}
\section{Fundamental Representations}
\setcounter{equation}{0}
\renewcommand{\theequation}{\thesection.\arabic{equation}}

\subsection{Affine weights}
\label{aff-wt}

Let $\ge$ be an affine Lie algebra. 
The sets $\mathfrak t$, 
$\{\al_i\}_{i\in I}$ 
and $\{\al^\vee_i\}_{i\in I}$ be as in \ref{KM}. 
We take ${\rm dim}\mathfrak t=\sharp I+1$.
Let $\del\in Q_+$ be the unique element 
satisfying $\{\lm\in Q|\lan \al^\vee_i,\lm\ran=0
\text{ for any }i\in I\}=\bbZ\del$
and ${\bf c}\in \ge$ be the canonical central element
satisfying $\{h\in Q^\vee|\lan h,\al_i\ran=0
\text{ for any }i\in I\}=\bbZ c$.
We write (\cite[6.1]{Kac})
\[
{\bf c}=\sum_i a_i^\vee \al^\vee_i,\qq
\del=\sum_i a_i\al_i.
\]
Let $(\q,\q)$ be the non-degenerate
$W$-invariant symmetric bilinear form on $\mathfrak t^*$
normalized by $(\del,\lm)=\lan {\bf c},\lm\ran$
for $\lm\in\frak t^*$.
Let us set $\tt^*_{\rm cl}:=\tt^*/\bbC\del$ and let
${\rm cl}:\tt^*\longrightarrow \tt^*_{\rm cl}$
be the canonical projection. 
Here we have 
$\tt^*_{\rm cl}\cong \oplus_i(\bbC \al^\vee_i)^*$.
Set $\tt^*_0:=\{\lm\in\tt^*|\lan {\bf c},\lm\ran=0\}$,
$(\tt^*_{\rm cl})_0:={\rm cl}(\tt^*_0)$. 
Since $(\del,\del)=0$, we have a positive-definite
symmetric form on $\tt^*_{\rm cl}$ 
induced by the one on 
$\tt^*$. 
Let $\Lm_i\in \tt^*_{\rm cl}$ $(i\in I)$ be a classical 
weight such that $\lan \al^\vee_i,\Lm_j\ran=\del_{i,j}$, which 
is called a fundamental weight.
We choose 
$P$ so that $P_{\rm cl}:={\rm cl}(P)$ coincides with 
$\oplus_{i\in I}\bbZ\Lm_i$ and 
we call $P_{\rm cl}$ a 
{\it classical weight lattice}.

\subsection{Fundamental representation 
$W(\varpi_1)$}
\label{fundamental}

Let $c=\sum_{i}a_i^\vee \al^\vee_i$ be the canonical
central element in an affine Lie algebra $\ge$
(see \cite[6.1]{Kac}), 
$\{\Lm_i|i\in I\}$ the set of fundamental 
weight as in the previous section
and $\varpi_1:=\Lm_1-a^\vee_1\Lm_0$ the
(level 0)fundamental weight.

Let $V(\varpi_1)$ be the extremal weight module
of $\uq$
associated with $\varpi_1$ (\cite{K0}) and 
$W(\varpi_1)\cong 
V(\varpi_1)/(z_1-1)V(\varpi_1)$ the 
fundamental representation of $\uqp$
where $z_1$ is a 
$\uqp$-linear automorphism on $V(\varpi_1)$
 (see \cite[Sect 5.]{K0}). 

By \cite[Theorem 5.17]{K0}, $W(\varpi_1)$ is an
finite-dimensional irreducible integrable 
$\uqp$-module and has a global basis
with a simple crystal. Thus, we can consider 
the specialization $q=1$ and obtain the 
finite-dimensional $\ge$-module $W(\varpi_1)$, 
which we call a fundamental representation
of $\ge$ and use the same notation as above.

We shall present the explicit form of 
$W(\varpi_1)$ for $\ge=\TY(D,3,4)$.
\subsection{$W(\varpi_1)$ for $\TY(D,3,4)$}
The Cartan matrix $A=(a_{i,j})_{i,j=0,1,2}$ of type 
$\TY(D,3,4)$ is as follows:
\[
 A=\begin{pmatrix}2&-1&0\\
-1&2&-3\\0&-1&2
\end{pmatrix}.
\]
Then the simple roots are 
\[
 \al_0=2\Lm_0-\Lm_1+\del,\q
\al_1=-\Lm_0+2\Lm_1-\Lm_2,\q
\al_2=-3\Lm_1+2\Lm_2, 
\]
and the Dynkin diagram is:
\[\SelectTips{cm}{}
\xymatrix{
*{\bigcirc}<3pt> \ar@{-}[r]_<{0} 
& *{\bigcirc}<3pt> \ar@3{<-}[r]_<{1}
& *{\bigcirc}<6pt>\ar@{}_<{\,\,\,\,\,\,2}
}
\]

The fundamental representation $W(\varpi_1)$ of type $\TY(D,3,4)$
is an 8-dimensional module with the basis,
\[
 \{v_1,v_2, v_3,v_0, {\emptyset}, v_{\ovl 3},
v_{\ovl 2},v_{\ovl 1}\}.
\]
The explicit form of $W(\varpi_1)$ is given 
in \cite{KMOY}.
\begin{eqnarray*}
&&\hspace{-30pt}{\rm wt}(v_1)=\Lm_1-2\Lm_0,\,\,
{\rm wt}(v_2)=-\Lm_0-\Lm_1+\Lm_2,\,\,
{\rm wt}(v_3)=-\Lm_0+2\Lm_1-\Lm_2,\\
&&\hspace{-30pt}{\rm wt}(v_{\ovl i})=
-{\rm wt}(v_i)\,\,(i=1,\cd,3),\,\,
{\rm wt}(v_0)={\rm wt}(\emptyset)=0.
\end{eqnarray*}
\def\bv#1{\fsquare(5mm,#1)}
The actions of $e_i$ and $f_i$ on these basis vectors
are given as follows:
\begin{eqnarray*}
&&\hspace{-30pt}
f_0\left(v_0,v_{\ovl 3}, v_{\ovl 2}, v_{\ovl 1},\emptyset\right)
=\left(v_1,v_2,v_3,\emptyset +\frac{1}{2}v_0,\frac{3}{2}v_1
\right),\\
&&\hspace{-30pt}f_1\left(v_1,v_3,v_0,v_{\ovl 2}, 
\emptyset\right)
=\left(v_2,v_0,2v_{\ovl 3}, v_{\ovl 1}\right),\\
&&\hspace{-30pt}f_2\left(v_2,v_{\ovl 3}\right)=
\left(v_3,v_{\ovl 2}\right),\\
&&\hspace{-30pt}e_0\left(v_1,v_2,v_3,v_0,\emptyset\right)=
\left(\emptyset+\frac{1}{2}v_0,v_{\ovl 3}, v_{\ovl 2},
v_{\ovl 1},\frac{3}{2}v_{\ovl 1}\right),\\
&&\hspace{-30pt}e_1\left(v_2,v_0,v_{\ovl 3},v_{\ovl 1}
\right)=
\left(v_1,2v_3,v_0,v_{\ovl 2}\right),\\
&&\hspace{-30pt}e_2\left(v_3,v_{\ovl 2}\right)=
\left(v_2,v_{\ovl 3}\right),
\end{eqnarray*}
where we only give non-trivial actions 
and the other actions are trivial.
We can easily check that these define 
the module $W(\varpi_1)$ by direct calculations.
\renewcommand{\thesection}{\arabic{section}}
\section{Affine Geometric Crystal $\cV_1(\TY(D,3,4))$}
\setcounter{equation}{0}
\renewcommand{\theequation}{\thesection.\arabic{equation}}

We shall construct the 
affine geometric crystal $\cV(\TY(D,3,4))$ in $W(\varpi_1)$
explicitly.

For $\xi\in (\frt^*_{\rm cl})_0$, let $t(\xi)$ be the 
translation as in \cite[Sect 4]{K0}.
Then we have 
\begin{eqnarray*}
&& t(\wtil\varpi_1)=s_0s_1s_2s_1s_2s_1=:w_1,\\
&& t(\text{wt}(v_{\ovl 2}))=s_2s_1s_2s_1s_0s_1=:w_2,
\end{eqnarray*}
Associated with these Weyl group elements $w_1$ and $w_2$,
we define subsets $\cV_1=\cV_1(\TY(D,3,4))$ and 
$\cV_2=\cV_2(\TY(D,3,4))\subset W(\varpi_1)$ 
respectively:
\begin{eqnarray*}
&&\hspace{-30pt}\cV_1:=\{V_1(x)
:=Y_0(x_0)Y_1(x_1)Y_2(x_2)Y_1(x_3)Y_2(x_4)Y_1(x_5)
v_1\,\,\vert\,\,x_i\in\bbC^\times,(0\leq i\leq 5)\},\\
&&\hspace{-30pt}\cV_2:=\{V_2(y):=
Y_2(y_2)Y_1(y_1)Y_2(y_4)Y_1(y_3)Y_0(y_0)Y_1(y_5)
v_{\ovl 2}\,\,\vert\,\,y_i\in\bbC^\times,(0\leq i\leq 5)\}.
\end{eqnarray*}
The set $\cV_1$ (resp. $\cV_2$) has a natural 
$G_2$ (resp. $A_2$)-geometric crystal structure
(see\cite{KNO})
for $i=1,2$ (resp. $i=0,1$), which is denoted by
$\Phi_1:=(\cV_1,\{e_1,e_2\},\{\gamma_1,\gamma_2\},
\{\vep_1,\vep_2\})$
(resp. $\Phi_2:=
(\cV_2,\{\ovl e_0,\ovl e_1\},
\{\ovl\gamma_0,\ovl\gamma_1\},
\{\ovl\vep_0,\ovl\vep_1\})$).
By the explicit forms of $f_i$'s on $W(\varpi_1)$
as above, we have $f_0^3=0$, $f_1^3=0$ and $f_2^2=0$ 
and then 
\begin{equation}
Y_i(c)=(1+\frac{f_i}{c}+\frac{f_i^2}{2c^2})\al_i^\vee(c)
\,\,(i=0,1),\q
Y_2(c)=(1+\frac{f_2}{c})\al_2^\vee(c).
\end{equation}
Thus, we can get explicit forms of $V_1(x)\in\cV_1$ 
and $V_2(y)\in\cV_2$. 
Set 
\begin{eqnarray*}
&&V_1(x)=\sum_{1\leq i\leq 3}\left(X_iv_{i}+X_{\ovl i}
v_{\ovl i}\right)+X_{0}v_0+X_\emptyset
\emptyset,\\
&&V_2(y)=\sum_{1\leq i\leq 3}\left(Y_iv_{i}+Y_{\ovl i}
v_{\ovl i}\right)+Y_{0}v_0+Y_\emptyset
\emptyset.
\end{eqnarray*}
Then by direct calculations, we have 
\begin{lem}
\label{XY}
The rational function $X_1,X_2,\cd,$ and $Y_1,Y_2,\cd$
are given as:
\begin{eqnarray*}
&&\hspace{-20pt}
X_1=
1 + \frac{{x_1}\,{x_3}\,{x_5}}{{{x_0}}^2} + 
  \frac{{x_3} + {x_1}\,
\left( \frac{{{x_3}}^2}{{x_2}} + 
\frac{{x_4}}{{x_3}} 
+ {x_5} \right) }{{x_0}}
,\\
&&\hspace{-20pt}X_2=
\frac{{x_2}}{{{x_1}}^2} + \frac{2\,{x_3}}{{x_1}} + 
  \frac{{{x_3}}^2}{{x_2}} + \frac{{x_4}}{{x_3}} + {x_5} + 
  \frac{{x_3}\,{x_5} 
+ \frac{{x_2}\,
\left( \frac{{x_4}}{{x_3}} 
+ {x_5} \right) }{{x_1}}}{{x_0}},
\\
&&\hspace{-20pt}X_3=
{x_1} + \frac{{{x_1}}^2\,\left( \frac{{{x_3}}^2}{{x_2}} + 
\frac{{x_4}}{{x_3}} 
+ {x_5} \right) }{{x_0}},\qq
X_{0}=
\frac{{x_0}}{2} + {x_3} 
+ {x_1}\,
\left( \frac{{{x_3}}^2}{{x_2}} 
+ \frac{{x_4}}{{x_3}} 
+ {x_5} \right),\\
&&\hspace{-20pt}
X_{\ovl 3}={x_0}\,\left( \frac{{x_2}}{{{x_1}}^2} + \frac{2\,{x_3}}{{x_1}} + 
\frac{{{x_3}}^2}{{x_2}} 
+ \frac{{x_4}}{{x_3}} + {x_5} \right),\q
X_{\ovl 2}=x_0x_1,\q
X_{\ovl 1}=x_0^2,\q
X_\emptyset=x_0,
\end{eqnarray*}
\begin{eqnarray*}
&&\hspace{-20pt}Y_1=y_1y_3,
\q Y_2=
{y_2}\,\left( {y_3} + \frac{{y_4}}{{y_1}} \right),\q
Y_3=
{y_3} + \frac{{y_4}}{{y_1}} + 
  \frac{{{y_1}}^2\,\left( 1 + \frac{{{y_3}}^2\,{y_5}}{{y_4}} \right) }
   {{y_2}}, \\
&&\hspace{-20pt}Y_0=
\frac{{y_0}}{2} + {y_3}\,{y_5} + 
  {y_1}\,
\left( 1 + \frac{{{y_3}}^2\,{y_5}}{{y_4}} \right), \\
&&\hspace{-20pt}Y_{\ovl 3}=
{y_2}\,\left( 1 + \frac{{{y_3}}^2\,{y_5}}{{y_4}} + 
    \frac{{y_4}\,\left( \frac{{y_0}}{{y_3}} + {y_5} \right) }
     {{{y_1}}^2} + \frac{2\,\left( \frac{{y_0}}{2} + {y_3}\,{y_5} \right) }
     {{y_1}} \right),\\
&&\hspace{-20pt}Y_{\ovl 2}=
1 + \frac{{{y_3}}^2\,{y_5}}{{y_4}} + 
  \frac{{y_4}\,\left( \frac{{y_0}}{{y_3}} + {y_5} \right) }{{{y_1}}^2} + 
  \frac{2\,\left( \frac{{y_0}}{2} + {y_3}\,{y_5} \right) }{{y_1}} + \frac{{y_1}\,\left( \frac{{y_0}}{{y_3}} + {y_5} +
\frac{{y_0}\,{y_3}\,{y_5}}{{y_4}} \right) }{{y_2}},\\
&&\hspace{-20pt}Y_{\ovl 1}=
\frac{{y_0}}{{y_3}} + {y_5} + 
  \frac{{y_0}\,{y_3}\,{y_5}}{{y_4}} + 
  \frac{\frac{{{y_0}}^2}{{y_3}} + {y_0}\,{y_5}}{{y_1}},\q
Y_\emptyset=y_0.
\end{eqnarray*}
\end{lem}
Now for a given $x=(x_0,\cd,x_5)$ we solve the equation 
\begin{equation}
 V_2(y)=a(x)V_1(x),
\label{eq}
\end{equation}
where $a(x)$ is a rational function in $x=(x_0,\cd,x_5)$.
Though this equation is over-determined, 
it can be solved uniquely and the explicit form of 
solution is as follows:
\begin{pro}
We have the rational function $a(x)$ and 
the unique solution of (\ref{eq}):
\begin{equation}
\begin{split}
&a(x)=\frac{1}{x_0\,{x_1}\,{x_2}\,
    {{x_3}}^2\,{x_4}}({x_0}\,{x_2}\,{x_3}\,{x_4} + 
    {x_2}\,{{x_3}}^2\,{x_4} + {x_1}\,{{x_3}}^3\,{x_4} + 
    {x_1}\,{x_2}\,{{x_4}}^2 \\& + 
    {x_0}\,{x_2}\,{{x_3}}^2\,{x_5} + 
    {x_0}\,{x_1}\,{{x_3}}^3\,{x_5}
     + 2\,{x_1}\,{x_2}\,{x_3}\,{x_4}\,{x_5} + 
    {x_1}\,{x_2}\,{{x_3}}^2\,{{x_5}}^2),\\
&y_0=a(x)x_0,\q
y_1=\frac{P}{{x_0}\,{x_1}\,{{x_2}}^2\,{{x_3}}^3\,{x_4}}
,\\
&
y_2=\frac{{x_2}}{{{x_1}}^3} + \frac{3\,{x_3}}{{{x_1}}^2} + 
  \frac{{{x_3}}^3}{{{x_2}}^2} + 
  \frac{{\left( {x_4} + {x_3}\,{x_5} \right) }^3}{{{x_3}}^3\,{x_4}} + 
  \frac{2\,{x_4} + 3\,{x_3}\,{x_5}}{{x_2}} + 
  \frac{3\,\left( {{x_3}}^3 + {x_2}\,{x_4} + 
       {x_2}\,{x_3}\,{x_5} \right) }{{x_1}\,{x_2}\,{x_3}},
\\
&
y_3=\frac{Q \,R }
    {x_0^2 P },\q
y_4=\frac{R^3}
    {S },\q
y_5=\frac{{x_0}\,{x_5}\,R}{Q },
\end{split}
\label{x->y}
\end{equation}
where 
\begin{eqnarray*}
P&=&{{x_1}}^2\,{{x_3}}^6\,{x_4} + 
    {x_1}\,{x_2}\,{{x_3}}^3\,{x_4}\,
     \left( 2\,{{x_3}}^2 + 2\,{x_1}\,{x_4} + 
       3\,{x_1}\,{x_3}\,{x_5} \right) 
 + {{x_2}}^2\,( {{x_3}}^4\,{x_4} + {{x_1}}^2\,{{x_4}}^3 \\
&&\hspace{-30pt}
+ 
 3\,{{x_1}}^2\,{x_3}\,{{x_4}}^2\,{x_5} + 
       {x_1}\,{{x_3}}^3\,{x_5}\,
        \left( 2\,{x_4} + {x_1}\,{{x_5}}^2 \right)  + 
       {x_1}\,{{x_3}}^2\,{x_4}\,
        \left( 2\,{x_4} + 3\,{x_1}\,{{x_5}}^2 \right))\\
&&\hspace{-30pt}  + 
    {x_0}\,{x_2}\,{x_3}\,
     \left( {x_1}\,{{x_3}}^3\,{x_4} + 
       {x_2}\,\left( {x_1}\,{{x_4}}^2 + 
          2\,{x_1}\,{x_3}\,{x_4}\,{x_5} + 
 {{x_3}}^2\,\left( {x_4} + {x_1}\,{{x_5}}^2 \right)  \right)  \right)
,\\
Q&=& {{x_0}}^2\,{x_2}\,{x_3} + 
      {x_1}\,{x_2}\,{{x_3}}^2\,{x_5} + 
      {x_0}\,\left( {x_1}\,{{x_3}}^3 + 
         {x_2}\,\left( {{x_3}}^2 + {x_1}\,{x_4} + 
            {x_1}\,{x_3}\,{x_5} \right)  \right),\\
R&=&{x_1}\,{{x_3}}^3\,{x_4} + 
       {x_0}\,{x_3}\,\left( {x_1}\,{{x_3}}^2\,{x_5} + 
          {x_2}\,\left( {x_4} + {x_3}\,{x_5} \right)  \right) \\&& + 
       {x_2}\,\left( {x_1}\,{{x_4}}^2 + 
          2\,{x_1}\,{x_3}\,{x_4}\,{x_5} + 
 {{x_3}}^2\,\left( {x_4} + {x_1}\,{{x_5}}^2 \right)  \right),\\
S&=&{{x_0}}^3\,{{x_3}}^3\,{x_4}\,
    ( {{x_2}}^3\,{{x_3}}^3\,{x_4} + 
      {{x_1}}^3\,{{x_3}}^6\,{x_4} + 
      {{x_1}}^2\,{x_2}\,{{x_3}}^3\,{x_4}\,
       \left( 3\,{{x_3}}^2 + 2\,{x_1}\,{x_4} + 
         3\,{x_1}\,{x_3}\,{x_5} \right) \\
&&\hspace{-30pt} + 
      {x_1}\,{{x_2}}^2\,( 3\,{{x_3}}^4\,{x_4} + 
         {{x_1}}^2\,{{x_4}}^3 + 
         3\,{{x_1}}^2\,{x_3}\,{{x_4}}^2\,{x_5} + 
         3\,{x_1}\,{{x_3}}^2\,{x_4}\,
          \left( {x_4} + {x_1}\,{{x_5}}^2 \right) \\
&&\hspace{-30pt}
 +  {x_1}\,{{x_3}}^3\,{x_5}\,
  \left( 3\,{x_4} + {x_1}\,{{x_5}}^2 \right)  ))
\end{eqnarray*}
Furthermore, 
the morphism given by (\ref{x->y})
\begin{eqnarray*}
\ovl\sigma:&\cV_1\longrightarrow &\cV_2,\\
&V_1(x_0,\cd,x_5)\mapsto &V_2(y_0,\cd,y_5).
\end{eqnarray*}
is a bi-positive birational isomorphism, that is, 
there exists an inverse positive rational 
map $\ovl\sigma^{-1}:\cV_2\to\cV_1$ 
$(V_2(y)\mapsto V_1(x))$:
\begin{eqnarray*}
&&x_0=\frac{1}{{y_3}} + \frac{{y_5}}{{y_0}} + 
  \frac{{y_3}\,{y_5}}{{y_4}} + 
  \frac{\frac{{y_0}}{{y_3}} + {y_5}}{{y_1}},\\
&&x_1=\frac{1 + \frac{{{y_3}}^2\,{y_5}}{{y_4}} + 
 \frac{{y_4}\,\left( \frac{{y_0}}{{y_3}} + {y_5} \right) }
 {{{y_1}}^2} + \frac{{y_0} + 2\,{y_3}\,{y_5}}{{y_1}} + 
 \frac{{y_1}\,\left( \frac{{y_0}}{{y_3}} + {y_5} + 
 \frac{{y_0}\,{y_3}\,{y_5}}{{y_4}} \right) }
{{y_2}}}{{y_0}},\\
&&x_2=\frac{T}{{{y_0}}^3\,
    {{y_1}}^6\,{y_2}\,{{y_3}}^3\,{{y_4}}^3},\q
x_3=\frac{U\,V}{W},\q
x_4=\frac{{{y_1}}^3\,{{y_2}}^2\,{y_4}\,  {U }^3}{y_0^3T},\q
x_5=\frac{{{y_1}}^2\,{y_2}\,{y_4}\,{y_5}\,U }
{{y_0}\,V},
\end{eqnarray*}
where
\begin{eqnarray*}
&&\hspace{-30pt}
T={{y_0}}^3\,{{y_4}}^3\,( {{y_1}}^6\,{y_4} + 
       3\,{{y_1}}^4\,{y_2}\,{y_3}\,{y_4} + 
       3\,{{y_1}}^2\,{{y_2}}^2\,{{y_3}}^2\,{y_4} + 
       3\,{y_1}\,{{y_2}}^2\,{y_3}\,{{y_4}}^2 + 
     {{y_2}}^2\,{{y_4}}^3 \\
&&\hspace{-20pt}+ {{y_1}}^3\,{y_2}\,( {y_2}\,
{{y_3}}^3 + 2\,{{y_4}}^2 ) ) 
 + 3\,{{y_0}}^2\,{y_3}\,{{y_4}}^2\,
 ( {{y_1}}^2\,{y_2}\,{{y_3}}^2 
+ {{y_1}}^3\,{y_4} + 2\,{y_1}\,{y_2}\,{y_3}\,{y_4} 
+ {y_2}\,{{y_4}}^2 ) \\
&&\hspace{-20pt}\times ( {{y_1}}^3\,{y_4}\,{y_5} + 
       2\,{y_1}\,{y_2}\,{y_3}\,{y_4}\,{y_5}
 +   {y_2}\,{{y_4}}^2\,{y_5} + 
  {{y_1}}^2\,{y_2}\,( {y_4} + {{y_3}}^2\,{y_5} )  ) \\
&&\hspace{-20pt}+ {{y_3}}^3\,( 6\,{y_1}\,
{{y_2}}^2\,{y_3}\,{{y_4}}^5\,
        {{y_5}}^3 + {{y_2}}^2\,{{y_4}}^6\,{{y_5}}^3 
+ 3\,{{y_1}}^2\,{{y_2}}^2\,{{y_4}}^4\,{{y_5}}^2\,
        ( {y_4} + 5\,{{y_3}}^2\,{y_5} ) \\
&&\hspace{-20pt}
 + 2\,{{y_1}}^3\,{y_2}\,{{y_4}}^3\,{{y_5}}^2\,
  ( 6\,{y_2}\,{y_3}\,{y_4} +   
10\,{y_2}\,{{y_3}}^3\,{y_5} + {{y_4}}^2\,{y_5} )\\
&&\hspace{-20pt}
+  3\,{{y_1}}^5\,{y_2}\,{y_4}\,{y_5}\,
        ( 2\,{y_2}\,{y_3}\,
    {( {y_4} + {{y_3}}^2\,{y_5} ) }^2 
+ {{y_4}}^2\,{y_5}\,( {y_4} + 2\,{{y_3}}^2\,{y_5} ))\\
&&\hspace{-20pt}
  + {{y_1}}^6\,( {{y_4}}^4\,{{y_5}}^3 + 
      {{y_2}}^2\,{( {y_4} + {{y_3}}^2\,{y_5} ) }^3
+ {y_2}\,{y_3}\,{{y_4}}^2\,{{y_5}}^2\,
           ( 3\,{y_4} + 2\,{{y_3}}^2\,{y_5} )  ) \\
&&\hspace{-20pt}
 +  3\,{{y_1}}^4\,{y_2}\,{{y_4}}^2\,{y_5}\,
        ( 2\,{y_3}\,{{y_4}}^2\,{{y_5}}^2 + 
   {y_2}\,( {{y_4}}^2 + 6\,{{y_3}}^2\,{y_4}\,{y_5} + 
             5\,{{y_3}}^4\,{{y_5}}^2 ))) \\
&&\hspace{-20pt}
 + 3\,{y_0}\,{{y_3}}^2\,{y_4}\,
  ( 5\,{y_1}\,{{y_2}}^2\,{y_3}\,{{y_4}}^4\,{{y_5}}^2 + 
       {{y_2}}^2\,{{y_4}}^5\,{{y_5}}^2 + 
       2\,{{y_1}}^2\,{{y_2}}^2\,{{y_4}}^3\,{y_5}\,
        ( {y_4} + 5\,{{y_3}}^2\,{y_5} )\\
&&\hspace{-20pt} 
+ 2\,{{y_1}}^3\,{y_2}\,{{y_4}}^2\,{y_5}\,
        ( 3\,{y_2}\,{y_3}\,{y_4}
 +  5\,{y_2}\,{{y_3}}^3\,{y_5} + {{y_4}}^2\,{y_5} )  + 
       {{y_1}}^6\,{y_4}\,{y_5}\,
        ( {{y_4}}^2\,{y_5} \\
&&\hspace{-20pt}+ 
     {y_2}\,{y_3}\,( {y_4} + {{y_3}}^2\,{y_5} )  ) 
       + {{y_1}}^5\,{y_2}\,( {y_2}\,{y_3}\,
           {( {y_4} + {{y_3}}^2\,{y_5} ) }^2 
+ 2\,{{y_4}}^2\,{y_5}\,( {y_4} 
+ 2\,{{y_3}}^2\,{y_5} )  )\\
&&\hspace{-20pt}
  + {{y_1}}^4\,{y_2}\,{y_4}\,
        ( 5\,{y_3}\,{{y_4}}^2\,{{y_5}}^2+ 
    {y_2}\,( {{y_4}}^2 + 6\,{{y_3}}^2\,{y_4}\,{y_5} + 
             5\,{{y_3}}^4\,{{y_5}}^2 )  )  ),\\
&&\hspace{-30pt}
U=( {{y_0}}^2\,{y_4} + {y_1}\,{y_3}\,{y_4}\,{y_5} + 
      {y_0}\,( {y_3}\,{y_4}\,{y_5} + 
 {y_1}\,( {y_4} + {{y_3}}^2\,{y_5} 
 )  )  ),\\
&&\hspace{-30pt}
V= ( 3\,{y_1}\,{y_2}\,{{y_3}}^2\,{{y_4}}^2\,{{y_5}}^2 + 
      {y_2}\,{y_3}\,{{y_4}}^3\,{{y_5}}^2 + 
      {{y_1}}^2\,{y_2}\,{y_3}\,{y_4}\,{y_5}\,
       ( 2\,{y_4} + 3\,{{y_3}}^2\,{y_5} )\\
&&\hspace{-20pt}  + 
      {y_0}\,{y_4}\,( {{y_1}}^3\,{y_4}\,{y_5} + 
     2\,{y_1}\,{y_2}\,{y_3}\,{y_4}\,{y_5} + 
     {y_2}\,{{y_4}}^2\,{y_5} + 
     {{y_1}}^2\,{y_2}\,( {y_4} + {{y_3}}^2\,{y_5} ) 
      ) \\
&&\hspace{-20pt}
 + {{y_1}}^3\,( {y_3}\,{{y_4}}^2\,{{y_5}}^2 + 
 {y_2}\,{( {y_4} + {{y_3}}^2\,{y_5} ) }^2 )),\\
&&\hspace{-30pt}
W={{y_0}}^2\,{y_1}\,{{y_4}}^2\,
    ( {y_0}\,( {y_1}\,{y_2}\,{y_3}\,{y_4} + 
   {y_2}\,{{y_4}}^2 + {{y_1}}^3\,
   ( {y_4} + {{y_3}}^2\,{y_5} ))  + 
   {y_3}\,( {{y_1}}^3\,{y_4}\,{y_5}\\
&&\hspace{-20pt} 
+ 2\,{y_1}\,{y_2}\,{y_3}\,{y_4}\,{y_5} + 
   {y_2}\,{{y_4}}^2\,{y_5} + 
   {{y_1}}^2\,{y_2}\,( {y_4} + {{y_3}}^2\,{y_5}))).
\end{eqnarray*}
\end{pro}
\begin{proof}
By the direct calculations, we obtain the results. 
Indeed, certain computer softwares are useful to the calculations.
\end{proof}
Here we obtain the positive birational isomorphism 
$\ovl\sigma:\cV_1\longrightarrow \cV_2$ ($V_1(x)\mapsto V_2(y)$)
 and its inverse $\ovl\sigma^{-1}$ as above.
The actions of $\ovl e_0^c$ on $V_2(y)$ 
(respectively $\ovl \gamma_0(V_2(y))$ and 
$\ovl \vep_0(V_2(y)))$ are 
induced from the ones on 
$Y_2(y_2)Y_1(y_1)Y_2(y_4)Y_1(y_3)Y_0(y_0)Y_1(y_5)$ 
as an element of the geometric crystal $\cV_2$ 
since
$e_0v_{\ovl 2}=e_1v_{\ovl 2}=0$. 
Now, we define the action $e_0^c$ on $V_1(x)$ by
\begin{equation}
e_0^cV_1(x)=\ovl\sigma^{-1}\circ \ovl e_0^c\circ
\ovl\sigma(V_1(x))).
\label{e0}
\end{equation}
We also define $\gamma_0(V_1(x))$ 
and $\vep_0(V_1(x))$ by 
\begin{equation}
\gamma_0(V_1(x))=\ovl\gamma_0(\ovl\sigma(V_1(x))),\qq
\vep_0(V_1(x)):=\ovl\vep_0(\ovl\sigma(V_1(x))).
\label{wt0}
\end{equation}
\begin{thm}
Together with 
(\ref{e0}), (\ref{wt0}) on $\cV_1$, we obtain a
positive affine geometric crystal $\chi:=
(\cV_1,\{e_i\}_{i\in I},
\{\gamma_i\}_{i\in I},\{\vep_i\}_{i\in I})$
$(I=\{0,1,2\})$, whose explicit form is as follows:
first we have $e_i^c$, $\gamma_i$ and $\vep_i$
for $i=1,2$ from the formula (\ref{eici}), (\ref{vep-i})
and (\ref{gamma-i}).
\begin{eqnarray*}
&&\hspace{-30pt}
e_1^c(V_1(x))=V_1(x_0,\cC_1x_1,x_2,\cC_3x_3,x_4,\cC_5x_5),\,
e_2^c(V_1(x))=V_1(x_0,x_1,\cC_2x_2,x_3,\cC_4x_4,x_5),\\
&&\text{where}\\
&&\cC_1=\frac{\frac{c\,{x_0}}{{x_1}} 
+ \frac{{x_0}\,{{x_2}}}{{{x_1}}^2\,{x_3}} 
+\frac{{x_0}\,{{x_2}}\,{{x_4}}}{{{x_1}}^2\,
{{x_3}}^2\,{x_5}}}{\frac{{x_0}}{{x_1}} 
+ \frac{{x_0}\,{{x_2}}}{{{x_1}}^2\,{x_3}} + 
\frac{{x_0}\,{{x_2}}\,
{{x_4}}}{{{x_1}}^2\,{{x_3}}^2\,{x_5}}},\q
\cC_3=\frac{\frac{c\,{x_0}}{{x_1}} 
+ \frac{c\,{x_0}\,{{x_2}}}{{{x_1}}^2\,{x_3}} + 
\frac{{x_0}\,{{x_2}}\,{{x_4}}}{{{x_1}}^2
\,{{x_3}}^2\,{x_5}}}{
\frac{c\,{x_0}}{{x_1}} 
+ \frac{{x_0}\,{{x_2}}}{{{x_1}}^2\,{x_3}} + 
\frac{{x_0}\,{{x_2}}\,{{x_4}}}
{{{x_1}}^2\,{{x_3}}^2\,{x_5}}},\\
&&\cC_5=\frac{c\,\left( \frac{{x_0}}{{x_1}} + 
      \frac{{x_0}\,{{x_2}}}{{{x_1}}^2\,{x_3}} + 
\frac{{x_0}\,{{x_2}}\,{{x_4}}}{{{x_1}}^2
\,{{x_3}}^2\,{x_5}} \right)}{\frac{c\,{x_0}}{{x_1}} 
+ \frac{c\,{x_0}\,{{x_2}}}{{{x_1}}^2\,{x_3}} + 
\frac{{x_0}\,{{x_2}}\,{{x_4}}}
{{{x_1}}^2\,{{x_3}}^2\,{x_5}}},\,
\cC_2=\frac{\frac{c\,{x_1}^3}{{x_2}} 
+ \frac{{x_1}^3\,{x_3}^3}{{{x_2}}^2\,{x_4}}}
{\frac{{x_1}^3}{{x_2}} + \frac{{x_1}^3\,{x_3}^3}{{{x_2}}^2\,{x_4}}},\,
\cC_4=\frac{c\,\left( \frac{{x_1}^3}{{x_2}} + 
\frac{{x_1}^3\,{x_3}^3}{{{x_2}}^2\,{x_4}} \right) }{\frac{c\,{x_1}^3}
{{x_2}} + \frac{{x_1}^3\,{x_3}^3}{{{x_2}}^2\,{x_4}}},\\
&&\vep_1(V_1(x))={\frac{{x_0}}{{x_1}} 
+ \frac{{x_0}\,{{x_2}}}{{{x_1}}^2\,{x_3}} + 
\frac{{x_0}\,{{x_2}}\,
{{x_4}}}{{{x_1}}^2\,{{x_3}}^2\,{x_5}}},\q
\vep_2(V_1(x))={\frac{{x_1}^3}{{x_2}} 
+ \frac{{x_1}^3\,{x_3}^3}{{{x_2}}^2\,{x_4}}},\\
&&\gamma_1(V_1(x))=\frac{x_1^2x_3^2x_5^2}{x_0x_2x_4},\q
\gamma_2(V_1(x))=\frac{x_2^2x_4^2}{x_1^3x_3^3x_5^3}.
\end{eqnarray*}
We also have $e_0^c$, $\vep_0$ and $\gamma_0$ on $V_1(x)$
as:
\begin{eqnarray*}
&&e_0^c(V_1(x))=V_1(\frac{D}{c\cdot E}x_0,\frac{F}{c\cdot E}x_1,
\frac{G}{c^3\cdot E^3}x_2,\frac{D\cdot H}{c^2\cdot E\cdot F}x_3,
\frac{D^3}{c^3\cdot G}x_4,\frac{D}{c\cdot H}x_5),\\
&&\vep_0(V_1(x))=\frac{E}{{{x_0}}^3\,{{x_2}}\,{x_3}},\qq
\gamma_0(V_1(x))=\frac{x_0^2}{x_1x_3x_5},
\end{eqnarray*}
{where}
\begin{eqnarray*}
&&\hspace{-20pt}
D=c^2\,{{x_0}}^2\,{x_2}\,{x_3} + 
 {x_1}\,{x_2}\,{{x_3}}^2\,{x_5} + 
  c\,{x_0}\,\left( {x_1}\,{{x_3}}^3 + 
 {x_2}\,\left( {{x_3}}^2 + {x_1}\,{x_4} + 
 {x_1}\,{x_3}\,{x_5} \right)  \right)
,\\
&&\hspace{-20pt}
E={{x_0}}^2\,{x_2}\,{x_3} + 
 {x_1}\,{x_2}\,{{x_3}}^2\,{x_5} + 
 {x_0}\,\left( {x_1}\,{{x_3}}^3 + 
 {x_2}\,\left( {{x_3}}^2 + {x_1}\,{x_4} + 
 {x_1}\,{x_3}\,{x_5} \right)  \right),\\
&&\hspace{-20pt}
F= {x_2}\,{{x_3}}^2\,
 \left( {x_0} + {x_1}\,{x_5} \right)  + 
 c\,{x_0}\,\left( {x_0}\,{x_2}\,{x_3} + 
 {x_1}\,\left( {{x_3}}^3 + {x_2}\,{x_4} + 
 {x_2}\,{x_3}\,{x_5} \right)  \right),\\
&&\hspace{-20pt}
G= c^3\,{{x_0}}^6\,{{x_2}}^3\,{{x_3}}^3 + 
      3\,c^2\,{{x_0}}^5\,{{x_2}}^3\,{{x_3}}^4 + 
      3\,c^2\,{{x_0}}^5\,{x_1}\,{{x_2}}^2\,{{x_3}}^5 + 
      3\,c\,{{x_0}}^4\,{{x_2}}^3\,{{x_3}}^5\\
&& + 6\,c\,{{x_0}}^4\,{x_1}\,{{x_2}}^2\,{{x_3}}^6
 +  {{x_0}}^3\,{{x_2}}^3\,{{x_3}}^6 + 
      3\,c\,{{x_0}}^4\,{{x_1}}^2\,{x_2}\,{{x_3}}^7 + 
      3\,{{x_0}}^3\,{x_1}\,{{x_2}}^2\,{{x_3}}^7 \\
&&+  3\,{{x_0}}^3\,{{x_1}}^2\,{x_2}\,{{x_3}}^8
    +  {{x_0}}^3\,{{x_1}}^3\,{{x_3}}^9 + 
3\,c^3\,{{x_0}}^5\,{x_1}\,{{x_2}}^3\,{{x_3}}^2\,{x_4}
+ 6\,c^2\,{{x_0}}^4\,{x_1}\,
  {{x_2}}^3\,{{x_3}}^3\,{x_4}\\
&& + 
3\,c\,{{x_0}}^4\,{{x_1}}^2\,{{x_2}}^2\,{{x_3}}^4\,{x_4}
+3\,c^3\,{{x_0}}^4\,{{x_1}}^2\,
{{x_2}}^2\,{{x_3}}^4\,{x_4} 
+3\,c\,{{x_0}}^3\,{x_1}\,{{x_2}}^3\,{{x_3}}^4\,{x_4}\\
&& + 
3\,{{x_0}}^3\,{{x_1}}^2\,{{x_2}}^2\,{{x_3}}^5\,{x_4}
+3\,c^2\,{{x_0}}^3\,{{x_1}}^2\,
{{x_2}}^2\,{{x_3}}^5\,{x_4}
+ 2\,{{x_0}}^3\,{{x_1}}^3\,{x_2}\,{{x_3}}^6\,{x_4}\\
&& + 
  c^3\,{{x_0}}^3\,{{x_1}}^3\,{x_2}\,{{x_3}}^6\,{x_4}
+ 3\,c^3\,{{x_0}}^4\,{{x_1}}^2\,
{{x_2}}^3\,{x_3}\,{{x_4}}^2 + 3\,c^2\,{{x_0}}^3\,{{x_1}}^2\,{{x_2}}^3\,{{x_3}}^2\,{{x_4}}^2\\
&& + 
 {{x_0}}^3\,{{x_1}}^3\,{{x_2}}^2\,{{x_3}}^3\,{{x_4}}^2 + 
2\,c^3\,{{x_0}}^3\,{{x_1}}^3\,{{x_2}}^2\,
{{x_3}}^3\,{{x_4}}^2
 +  c^3\,{{x_0}}^3\,{{x_1}}^3\,{{x_2}}^3\,{{x_4}}^3\\
&&+3\,c^3\,{{x_0}}^5\,{x_1}\,{{x_2}}^3\,{{x_3}}^3\,{x_5} + 
   9\,c^2\,{{x_0}}^4\,{x_1}\,{{x_2}}^3\,{{x_3}}^4\,{x_5} 
 6\,c^2\,{{x_0}}^4\,{{x_1}}^2\,{{x_2}}^2\,
{{x_3}}^5\,{x_5}\\
&& + 9\,c\,{{x_0}}^3\,{x_1}\,{{x_2}}^3\,{{x_3}}^5\,{x_5} + 
12\,c\,{{x_0}}^3\,{{x_1}}^2\,{{x_2}}^2\,{{x_3}}^6\,{x_5} 
+  3\,{{x_0}}^2\,{x_1}\,{{x_2}}^3\,{{x_3}}^6\,{x_5}\\
&& + 3\,c\,{{x_0}}^3\,{{x_1}}^3\,{x_2}\,{{x_3}}^7\,{x_5} + 
  6\,{{x_0}}^2\,{{x_1}}^2\,{{x_2}}^2\,{{x_3}}^7\,{x_5} + 
      3\,{{x_0}}^2\,{{x_1}}^3\,{x_2}\,{{x_3}}^8\,{x_5} \\
&&+6\,c^3\,{{x_0}}^4\,{{x_1}}^2\,{{x_2}}^3\,
{{x_3}}^2\,{x_4}\, {x_5} 
+ 12\,c^2\,{{x_0}}^3\,{{x_1}}^2\,{{x_2}}^3\,{{x_3}}^3\,
  {x_4}\,{x_5} + 3\,c\,{{x_0}}^3\,{{x_1}}^3\,{{x_2}}^2\,
       {{x_3}}^4\,{x_4}\,{x_5} \\
&&+3\,c^3\,{{x_0}}^3\,{{x_1}}^3\,{{x_2}}^2\,
{{x_3}}^4\,{x_4}\,{x_5} 
+ 6\,c\,{{x_0}}^2\,{{x_1}}^2\,{{x_2}}^3\,{{x_3}}^4\,
{x_4}\,{x_5} + 3\,{{x_0}}^2\,{{x_1}}^3\,{{x_2}}^2\,
       {{x_3}}^5\,{x_4}\,{x_5}
\end{eqnarray*}
\begin{eqnarray*}
&& + 3\,c^2\,{{x_0}}^2\,{{x_1}}^3\,
{{x_2}}^2\,{{x_3}}^5\,{x_4}\,
       {x_5} + 3\,c^3\,{{x_0}}^3\,
{{x_1}}^3\,{{x_2}}^3\,{x_3}\,
{{x_4}}^2\,{x_5} + 3\,c^2\,{{x_0}}^2\,
{{x_1}}^3\,{{x_2}}^3\,
{{x_3}}^2\,{{x_4}}^2\,{x_5} \\
&&+  3\,c^3\,{{x_0}}^4\,{{x_1}}^2\,
{{x_2}}^3\,{{x_3}}^3\,{{x_5}}^2 + 
 9\,c^2\,{{x_0}}^3\,{{x_1}}^2\,{{x_2}}^3\,
{{x_3}}^4\,{{x_5}}^2 + 
 3\,c^2\,{{x_0}}^3\,{{x_1}}^3\,
{{x_2}}^2\,{{x_3}}^5\,{{x_5}}^2 \\
&&+ 9\,c\,{{x_0}}^2\,{{x_1}}^2\,{{x_2}}^3\,
{{x_3}}^5\,{{x_5}}^2 + 
 6\,c\,{{x_0}}^2\,{{x_1}}^3\,{{x_2}}^2\,{{x_3}}^6\,{{x_5}}^2 + 
3\,{x_0}\,{{x_1}}^2\,{{x_2}}^3\,{{x_3}}^6\,{{x_5}}^2 \\
&&+ 3\,{x_0}\,{{x_1}}^3\,{{x_2}}^2\,{{x_3}}^7\,{{x_5}}^2 
+3\,c^3\,{{x_0}}^3\,{{x_1}}^3\,{{x_2}}^3\,
{{x_3}}^2\,{x_4}\,{{x_5}}^2 + 
6\,c^2\,{{x_0}}^2\,{{x_1}}^3\,{{x_2}}^3\,{{x_3}}^3\,
   {x_4}\,{{x_5}}^2\\
&& + 3\,c\,{x_0}\,{{x_1}}^3\,{{x_2}}^3\,
       {{x_3}}^4\,{x_4}\,{{x_5}}^2 
+  c^3\,{{x_0}}^3\,{{x_1}}^3\,{{x_2}}^3\,{{x_3}}^3\,
{{x_5}}^3 +   3\,c^2\,{{x_0}}^2\,{{x_1}}^3\,
{{x_2}}^3\,{{x_3}}^4\,{{x_5}}^3 \\
&&+3\,c\,{x_0}\,{{x_1}}^3\,{{x_2}}^3\,
{{x_3}}^5\,{{x_5}}^3 + 
 {{x_1}}^3\,{{x_2}}^3\,{{x_3}}^6\,{{x_5}}^3,\\
&&\hspace{-20pt}
H= c\,{{x_0}}^2\,{x_2}\,{x_3} + 
 {x_0}\,{x_2}\,{{x_3}}^2 + {x_0}\,{x_1}\,{{x_3}}^3 + 
 {x_0}\,{x_1}\,{x_2}\,{x_4} + 
 c\,{x_0}\,{x_1}\,{x_2}\,{x_3}\,{x_5} + 
 {x_1}\,{x_2}\,{{x_3}}^2\,{x_5}.
\end{eqnarray*}
\end{thm}
\begin{proof}
Since the positivity is trivial from the above formula,
it suffices to show that $\chi:=
(\cV_1,\{e_i^c\}_{i\in I},
\{\gamma_i\}_{i\in I},$ $\{\vep_i\}_{i\in I})$ 
satisfies the relations in Definition \ref{def-gc}.
Indeed, $\cV_1$ is $G_2$-geometric crystal for $i=1,2$
and then we may see the cases related to $i=0$,
that is, 
\begin{eqnarray}
&&\gamma_0(e_i^c(x))=c^{a_{i,0}}\gamma_0(x),\q
\gamma_i(e_0^c(x))=c^{a_{0,i}}\gamma_i(x),\q
\vep_0(e_0^c(x))=c^{-1}\vep_0(x),\label{0g-ep}\\
&&e_0^ce_1^{cd}e_0^d=e_1^de_0^{cd}e_1^c,\q
e_0^ce_2^d=e_2^de_0^c.\label{012}
\end{eqnarray}
In order to show these relations, we need:
\begin{lem}
$\ovl\sigma\circ e_1^c(x)
=\ovl e_1^c\circ\ovl\sigma(x)$ and 
$\gamma_1(x)=\ovl\gamma_1\circ\osigma(x)$
for $x\in\cV_1$.
\end{lem}
This lemma is obtained by calculating directly.

The relations in (\ref{0g-ep}) for $i=2$ are shown:
\[
 \gamma_2(e_0^c(x)=\left(\frac{G}{c^3E^3}x_2
\cdot\frac{D^3}{c^3G}x_4\right)^2
\left(\frac{F}{cE}x_1\cdot\frac{DH}{c^2EF}x_3
\frac{D}{cH}x_5\right)^{-3}
=\frac{x_2^2x_4^2}{x_1^3x_3^3x_5^3}=\gamma_2(x),
\]
and $\gamma_0(e_2^c(x))=\gamma_0(x)$ since 
$\gamma_0(x)$  does not depend on $x_2,\,\,x_4$.
The other relations in (\ref{0g-ep}) 
are shown by using the lemma:
\begin{eqnarray*}
&& \gamma_0(e_1^c(x))
=\ovl\gamma_0(\osigma\circ e_1^c(x))
=\ovl\gamma_0(\ovl e_1^c\circ\osigma(x))
=c^{-1}\ovl\gamma_0(\osigma(x))=c^{-1}\gamma_0(x),\\
&&\gamma_0(e_0^c(x))
=\ovl\gamma_0(\osigma(\osigma^{-1}
\circ\ovl e_0^c\circ\osigma(x)))
=\ovl\gamma_0(\ovl e_0^c\circ\osigma(x))
=c^2\ovl\gamma_0(\osigma(x))=c^2\gamma_0(x),\\
&&\gamma_1(e_0^c(x))
=\ovl\gamma_1(\osigma(\osigma^{-1}
\circ\ovl e_0^c\circ\osigma(x)))
=\ovl\gamma_1(\ovl e_0^c\circ\osigma(x))
=c^{-1}\ovl\gamma_1(\osigma(x))=c^{-1}\gamma_1(x),\\
&&\vep_0(e_0^c(x))
=\ovl\vep_0(\osigma(\osigma^{-1}
\circ\ovl e_0^c\circ\osigma(x)))
=\ovl\vep_0(\ovl e_0^c\circ\osigma(x))
=c^{-1}\ovl\vep_0(\osigma(x))=c^{-1}\vep_0(x).
\end{eqnarray*}
In (\ref{012}) the first relation is shown by
using the lemma:
\[
 e_0^ce_1^{cd}e_0^d=(\osigma^{-1}\ovl e_0^c\osigma)
(\osigma^{-1}\ovl e_1^{cd}\osigma)
(\osigma^{-1}\ovl e_0^d\osigma)
=\osigma^{-1}\ovl e_0^c\ovl e_1^{cd}\ovl e_0^d\osigma
=\osigma^{-1}\ovl e_1^d\ovl e_0^{cd}\ovl e_1^c\osigma
=e_1^de_0^{cd}e_1^c.
\]
The last relation $e_0^ce_2^d=e_2^de_0^c$ is shown by
direct calculations. For the calculations,
some computer software is very useful.
Then we know that $\chi$ is an affine
geometric crystal of type $\TY(D,3,4)$.
\end{proof}
Here we denote the positive structure on $\chi$ by 
$\theta:\cV_1\longrightarrow T$. Then by Corollary \ref{cor-posi}
we obtain the ultra-discretization ${\mathcal UD}(\chi,T,\theta)$, 
which is a Kashiwara's crystal. In \cite{KNO},
\cite{N4},  we show that
such crystal is isomorphic to the limit of 
certain perfect crystal for 
the Langlands dual algebra. So we present the following 
conjecture:
\begin{conj}
\label{conjecture}
The crystal ${\mathcal UD}(\chi,T,\theta)$ as above 
is the limit of coherent family of 
perfect crystals of type $\TY(G,1,2)$ in \cite{Y}.
\end{conj}

\bibliographystyle{amsalpha}

\end{document}
